\documentclass{tranx}
\usepackage{natbib, amssymb,latexsym, amscd}
\usepackage[all]{xy}
\usepackage{graphicx}

 \newtheorem{theorem}{Theorem}[section]
 \newtheorem{cor}[theorem]{Corollary}
 \newtheorem{lemma}[theorem]{Lemma}
 \newtheorem{proposition}[theorem]{Proposition} \theoremstyle{definition}
 \newtheorem{definition}[theorem]{Definition}
 \theoremstyle{definition}
 \newtheorem{example}[theorem]{Example}
 \theoremstyle{remark}
 
 \numberwithin{equation}{section}

\newcommand{\ben}{\begin{equation}}
\newcommand{\een}{\end{equation}}

\newcommand{\integer}{\ensuremath{{\mathbb Z}}}

\newcommand{\real}{\ensuremath{{\mathbb R}}}
\newcommand{\complex}{\ensuremath{{\mathbb C}}}

\newcommand{\gr}{\mathfrak}

\begin{document}

\title[Orbifold Virtual Cohomology of the Symmetric Product]{Orbifold Virtual Cohomology of the Symmetric Product}

\author[David Riveros and Bernardo Uribe]{David Riveros and Bernardo Uribe}
\thanks{The first author was partially supported by the ``Proyecto Semilla"
from the Universidad de los Andes. The second author was partially
supported by the ``Fondo de apoyo a investigadores jovenes" from
Universidad de los Andes and from grant N. 120440520139 of
COLCIENCIAS}

\address{Departamento de Matem\'{a}ticas, Universidad de los Andes,
Carrera 1 N. 18A - 10, Bogot\'a, COLOMBIA}
 \email{ dav-rive@uniandes.edu.co \\
buribe@uniandes.edu.co }
\begin{abstract}
The virtual cohomology of an orbifold is a ring structure on the
cohomology of the inertia orbifold whose product is defined via
the pull-push formalism and the Euler class of the excess
intersection bundle. In this paper we calculate the virtual
cohomology of a large family of orbifolds, including the symmetric
product.
\end{abstract}
\maketitle

\section{Introduction}

It was noticed in \cite{LupercioUribeXicotencatl2} that the ring
structure defined in the homology of the loop space of the
symmetric product orbifold (see \cite{LupercioUribeXicotencatl})
induces a ring structure on the cohomology of the inertia
orbifold, by restricting the structure to the constant loops. This
led the authors of \cite{LupercioUribeXicotencatl2} to define a
ring structure on the inertia orbifold of any orbifold that the
authors coined {\it{virtual cohomology}}. This cohomology is
defined via the pull-push formalism in as much as the same way
that the Chen-Ruan product for orbifolds is defined (see
\cite{ChenRuan, FantechiGottsche}). In
\cite{GonzalezLupercioSegoviaUribeXicotencatl} the relation
between the virtual  and the Chen-Ruan cohomology was clarified,
namely, that for an almost complex orbifold, its virtual
cohomology is isomorphic to the Chen-Ruan cohomology of its
cotangent orbifold.

In this paper we give an algorithm to calculate the virtual
cohomology of a large family of orbifolds. We first show that for
any global quotient orbifold $[Y/G]$, the virtual cohomology
$H^*_{virt}(Y,G; \integer)$ maps to the group ring $H^*(Y;
\integer)[G]$, and therefore, when this map is injective we can
see the virtual cohomology as a subring of the group ring. This
for example is the case when the inclusions of the fixed point
sets $Y^g \to Y$ induce a monomorphism in homology. We calculate
the virtual cohomology of these orbifolds by describing a set of
generators in the group ring. In the case of the symmetric product
we reduce the set of generators to the lower degree cohomology
classes of the fixed point sets of the transpositions.

This paper was motivated by the master's thesis of the first
author \cite{Riveros} where the virtual cohomology of the
symmetric product of spheres was calculated.  Lastly, we would
like to thank A. Cardona, E. Lupercio and M. Xicotencatl for
useful conversations and the anonymous referee for pointing out
some redundant relations on the presentation of the last example.

\section{Virtual Cohomology}
Let $[Y/G]$ be an orbifold with $Y$ differentiable, compact,
oriented and closed, and $G$ a finite group acting smoothly on $Y$
preserving the orientation. The inertia orbifold $I[M/G]$ is
defined as the orbifold
$$I[Y/G] = \left[ \left(\bigsqcup_{g \in G} Y^g \times \{g\}
\right) /G \right]$$ where $Y^g$ denotes the fixed point set of
the element $g$, we label the components with the elements of the group  and $G$ acts in the following way:
\begin{eqnarray} \nonumber
\left( \bigsqcup_{g \in G} Y^g  \times \{g\} \right)\times G & \to
& \bigsqcup_{g
\in G} Y^g \times \{g\} \\
((x,g),h) & \mapsto & (xh, h^{-1}gh). \label{G action}
\end{eqnarray}

From \cite{LupercioUribeXicotencatl2} we know that the virtual
intersection product defines a ring structure on the cohomology of
the inertia orbifold $I[Y/G]$, this ring is what the authors in \cite{LupercioUribeXicotencatl2} have called
{\it{virtual cohomology}}. Let's recall its definition.

Consider the groups
$$H^*(Y, G; \integer) := \bigoplus_{g \in G} H^*(Y^g; \integer) \times \{g\}$$
and for $g,h \in G$ define the maps
\begin{eqnarray*}
\times : H^*(Y^g; \integer) \times H^*(Y^h; \integer) & \to & H^{*}(Y^{gh};\integer)\\
(\alpha, \beta) & \mapsto & i_{gh!}(i_g^* \alpha \cdot i_h^*\beta
\cdot e(Y, Y^g, Y^h))
\end{eqnarray*}
where $i_g : Y^g\cap Y^h \to Y^g$, $i_h : Y^g\cap Y^h \to Y^h$ and
$i_{gh} : Y^g\cap Y^h \to Y^{gh}$ are the inclusion maps, $e(Y,
Y^g, Y^h)$ is the Euler class of the excess bundle of the
inclusions $Y^g \to Y \leftarrow Y^h$ (see \cite{Quillen}) and
$i_{gh!}$ is the pushforward map in cohomology.

 In
\cite{LupercioUribeXicotencatl2} it was required that the orbifold
 be almost complex with a compatible $G$ action, but for the
product to be well defined it is only necessary that the Euler
classes of the excess bundles be of even degree. This can be
achieved if for all $g_i \in G$ the fixed point sets
$$Y^{g_1,\dots,g_n}:=Y^{g_1}\cap \cdots \cap Y^{g_n} $$ are of even
dimension.

The group $G$ acts on $H^*(Y,G; \integer)$ in the following way:
for $g,h \in G$ and $\alpha \in H^*(Y^g ; \integer)$ we have
$$(\alpha, g)\cdot h := ((h^{-1})^*\alpha, h^{-1}gh).$$

\begin{definition} Let $[Y/G]$ be an orbifold  such that for all $g_i \in G$ the fixed point sets
$Y^{g_1,\dots,g_n} $ are even dimensional. Then, the group
$H^*(Y,G; \integer)$ together with the ring structure
\begin{eqnarray*}
\bullet : H^*(Y,G; \integer) \times H^*(Y,G; \integer) & \to & H^*(Y,G; \integer) \\
((\alpha,g) , (\beta, h)) & \mapsto & (\alpha \times \beta, gh)
\end{eqnarray*}
is what is called the {\it{virtual cohomology}} of the pair
$(Y,G)$; we will denote it by $H^*_{virt}(Y,G; \integer)$ .
Moreover, as the ring structure is $G$-equivariant with respect to
the action of $G$ on $H^*(Y,G; \integer)$, we define the virtual
cohomology of the orbifold $[Y/G]$ as the $G$ invariant part of
the ring $H^*_{virt}(Y,G; \real)$, i.e.
$$H^*_{virt}([Y/G]; \real) := H^*(Y,G; \real)^G.$$
\end{definition}

In what follows we will show how to calculate the virtual
cohomology for a large family of orbifolds.

For all $ g \in G$ let $f_g : Y^g \to Y$ be the inclusion of
manifolds and $$f_{g!}: H^*(Y^g; \integer) \to H^*(Y; \integer)$$
be the pushforward in cohomology. Consider the group ring $H^*(Y;
\integer)[G]$ of the group $G$ with coefficients in the ring
$H^*(Y; \integer)$, together with the $G$ action defined by
$$\left( \sum_i \alpha_i g_i \right)\cdot h := \sum_i ((h^{-1})^*
\alpha_i)h^{-1}g_ih.$$

\begin{theorem}
The inclusions $f_g : Y^g \to Y$ induce an equivariant  ring homomorphism from
the virtual cohomology to the group ring
\begin{eqnarray*}
f : H^*_{virt}(Y,G; \integer) & \to  & H^*(Y;\integer)[G] \\
(\alpha,g) & \mapsto & (f_{g!} \alpha)g.
\end{eqnarray*}
\end{theorem}
\begin{proof}
To show that the map $f$ is a ring homomorphism we only need to
check the commutativity of the following diagram
$$\xymatrix{H^*(Y^g; \integer) \times H^*(Y^h; \integer) \ar[rr]^{f_{g!} \times
f_{h!}} \ar[d]^{\times} && H^*(Y; \integer) \times H^*(Y; \integer) \ar[d]^\cdot \\
H^*(Y^{gh}; \integer) \ar[rr]_{f_{hg!}} && H^*(Y; \integer). }$$
Consider the diagram of inclusions
$$\xymatrix{Y & Y^g \ar[l]_{f_g} \\
Y^h \ar[u]^{f_h} & Y^{g,h} \ar[l]^{i_h} \ar[u]_{i_g} \ar[ul]_s.
}$$ It was proven in \cite[Lemma 16]{LupercioUribeXicotencatl2} by
an application of Quillen's excess intersection formula
\cite[Prop. 3.3]{Quillen} that for $\alpha \in H^*(Y^g; \integer)$
and $ \beta \in H^*(Y^h; \integer)$ one has
$$f_{g!}\alpha \cdot f_{h!} \beta = s_!( i_g^* \alpha \cdot i_h^*
\beta \cdot e(Y, Y^g, Y^h)).$$ Now, as $s = f_{gh}\circ i_{gh}$ we
have that $s_!= f_{gh!} \circ i_{gh!}$ and therefore
\begin{eqnarray*}
f_{g!}\alpha \cdot f_{h!} \beta & = & f_{gh!} (i_{gh!}( i_g^*
\alpha \cdot i_h^* \beta \cdot e(Y, Y^g, Y^h)))\\
& = & f_{gh!}( \alpha \times \beta).
\end{eqnarray*}
To check that the map $f$ is $G$-equivariant we simply consider
the inclusion
\begin{eqnarray*}
\psi: \bigsqcup_{g \in G} Y^g \times \{g\} & \to & \bigsqcup_{g
\in G}
Y^g \times \{g\}\\
(x,g ) & \mapsto & (f_gx, g).
\end{eqnarray*}
If we endow the space $\bigsqcup_{g \in G} Y^g \times \{g\}$ with
the same $G$-action as in \ref{G action} then the map $\psi$
becomes $G$-equivariant and we have the commutativity of the
following square:

$$\xymatrix{Y^g \times \{g\}  \ar[rr]^{f_{g}} \ar[d]^{h} && Y \times \{g\} \ar[d]^h \\
Y^{h^{-1}gh} \times \{h^{-1}gh\}  \ar[rr]_{f_{h^{-1}gh}} && Y
\times \{h^{-1}gh\}. }$$ Therefore we have that for $g,h \in G$
and $\alpha \in H^*(Y^g; \integer)$
\begin{eqnarray*}
(f_{g!} \alpha, g) \cdot h & = & \left( (h^{-1})^* f_{g!} \alpha,
h^{-1}gh \right) \\
& = & \left( f_{h^{-1}gh !} (h^{-1})^* \alpha, h^{-1}gh  \right)
\end{eqnarray*}
and this implies that the map $f$ is $G$-equivariant.

\end{proof}

\begin{cor}
If the inclusion maps in homology $f_{g*} : H_*(Y^g ; \integer)
\to H_*(Y; \integer)$ are injective for all $g \in G$, then the
map $$f: H^*_{virt}(Y,G; \integer) \to H^*(Y; \integer)[G]$$ is an
injective homomorphism of rings. Then the ring $H^*(Y,G;
\integer)$ can be calculated as the subring $f(H^*_{virt}(Y,G;
\integer))$ of $H^*(Y; \integer)[G]$.
\end{cor}

\begin{proof}
The pushforward $f_{g!} : H^*(Y^g; \integer) \to H^*(Y; \integer)$
in cohomology can be defined as the composition of the maps $PD
\circ f_{g*} \circ PD_g^{-1}$ where $PD : H_*(Y; \integer)
\stackrel{\cong}{\to} H^*(Y; \integer)$ and $PD_g : H_*(Y^g;
\integer) \stackrel{\cong}{\to} H^*(Y^g; \integer)$ are the
Poincar\'e duality isomorphisms. It follows that the maps $f_{g!}$
are injective.
\end{proof}
The above corollary will allow us to calculate the virtual
cohomology of a large family of orbifolds, as in the following
example.

\begin{example} \label{example projective space}
Consider the action of $\integer/p$ on the complex projective
space $\complex P^n$
\begin{eqnarray*}
\complex P^n \times \integer/p & \to & \complex P^n \\
([z_0: \dots : z_n], \lambda^i) & \mapsto & [z_0: \dots : z_{n-1}
: \lambda^i z_n]
\end{eqnarray*}
where the elements of $\integer/p$ are taken as $p$-th roots of
unity.
For $i\neq 0$  one has that the fixed point set of $\lambda^i$ is
$$(\complex P^n)^{\lambda^i} \cong \complex P^{n-1} \cup \{*\}$$ and therefore
if we only consider the connected component of $\complex P^{n-1}$ then the maps $f_{\lambda^i*}$ are all injective. The pushforward of
the inclusions are
\begin{eqnarray*}
f_{\lambda^i !} : H^*(\complex P^{n-1}; \integer) =
\integer[y]/\langle y^n \rangle & \to &
H^*(\complex P^{n}; \integer) = \integer[x]/ \langle x^{n+1}\rangle \\
y^j & \mapsto & x^{j+1},
\end{eqnarray*}
and
\begin{eqnarray*}
f_{\lambda^i !} : H^*(\{*\}; \integer) = \integer \langle z
\rangle  & \to &
H^*(\complex P^{n}; \integer) = \integer[x]/ \langle x^{n+1}\rangle \\
z & \mapsto & x^n.
\end{eqnarray*}

Therefore we have

$$H^*_{virt}(\complex P^n, \integer/p \ ; \integer) \cong \left( \integer[x]/ \langle x^{n+1} \rangle [1] \oplus
\bigoplus_{i=1}^{p-1} \left( x\integer[x]/ \langle x^{n+1}\rangle
\oplus \integer \langle z \rangle \right) [\lambda^i] \right),$$
where $zx=0$ and $z^2=0$.

 If we
take a closer look at the virtual cohomology generated by the inclusions of the $\complex P^{n-1}$'s, we can see that its
elements are truncated polynomials of maximum degree $n$, whose
coefficients are elements in the group ring $\integer[\integer/p]$
except for the constant term that it should be an integer, i.e.
$$ \{ P(x) \in
\integer[\integer/p][x]/\langle x^{n+1}\rangle | P(0) \in \integer \}$$ where
the ring structure is given by
multiplication of polynomials. If we add the classes coming from the inclusions of the points $*$ we obtain that
$$H^*_{virt}(\complex P^n, \integer/p \ ; \integer) \cong \{ P(x) \in
\integer[\integer/p][x]/\langle x^{n+1}\rangle | P(0) \in \integer
\} \oplus  \oplus_{i=1}^p \integer\langle z \rangle[\lambda^i]
/(xz, z^2).$$

Now, as the group $\integer/p$ is abelian and its action can be
factored through an action of $S^1$, we have that
$$H^*_{virt}(\complex P^n , \integer/p \ ; \real)^{\integer/p} =H^*_{virt}(\complex P^n ,
\integer/p \ ; \real).$$ Then
$$H^*_{virt}([\complex P^n/ \integer/p] ; \real) \cong \{ P(x) \in
\real[\integer/p][x]/\langle x^{n+1}\rangle | P(0) \in \real \}
 \oplus  \oplus_{i=1}^p \real\langle z \rangle[\lambda^i]
/(xz, z^2).$$
\end{example}

We have seen that for the case in which the homomorphisms $f_{g*}$
are injective, the virtual cohomology is isomorphic to the subring
$f(H^*_{virt}(Y,G; \integer))$ of the group ring $H^*(Y;
\integer)[G]$. In what follows we will find a set of generators
for $f(H^*_{virt}(Y,G; \integer))$.

Let $H^*(Y; \integer)[1_G]$ be the set of elements of the group
ring whose label is the identity $1_G$ of the group $G$.

\begin{proposition} \label{proposition generators G}
Suppose that all the homomorphisms $f_{g*}$ are injective and all
the $f_{g}^{*}$ are surjective. Denote by $1_g \in H^0(Y^g;
\integer)$ the identity of the ring $H^*(Y^g ; \integer)$. Then
the set
$$W:=H^*(Y; \integer)[1_G] \cup \{(f_{g!}1_g)g \ |\  g \in G \}$$
generates the ring $f(H^*_{virt}(Y,G; \integer))$.
\end{proposition}

\begin{proof}
It is clear that $W \subset f(H^*_{virt}(Y,G; \integer))$; we need
to prove that for any $a \in H^*(Y^g; \integer)$ the element
$(f_{g!}a) g$ can be generated with elements in $W$.

We know that the pullback $f_{g}^*$ is surjective. Therefore there
exists $b \in H^*(Y; \integer)$ such that $f_g^* b = a$. By the
module structure of the pushforward we have
$$f_{g!}(a) = f_{g!}(1_g a) = f_{g!}(1_g f_g^*b)= (f_{g!} 1_g)b$$
which implies that in the group ring
$$(1_g \ g)(b\  1_G)= \left( (f_{g!} 1_g)b\right)g  =(f_{g!}a) g.  $$
\end{proof}

\section{Symmetric Product}

It was shown in \cite{LupercioUribeXicotencatl2} that for an even
dimensional compact and closed manifold $M$, the virtual
cohomology of the orbifold $[M^n/ \gr{S}_n]$ is a subring of
 the string homology of the loop orbifold of the symmetric product
 (see \cite{LupercioUribeXicotencatl}).

In this section we will calculate the virtual cohomology of the
pair $(M^n, \gr{S}_n)$ in terms of the cohomology of $M$. As we
will make use of the Kunneth isomorphism we will restrict to real
coefficients. We would like to remark that if the manifold has
torsion free homology, all the calculations that follow can be
done with integer coefficients.

So, abusing the notation, we will talk indistinctly of the rings
$H^*(M^k ; \real)$ and $H^*(M ; \real)^{\otimes k}$.

We know that the diagonal inclusion $\Delta : M \to M\times M$
induces an injection $\Delta_* : H_*(M ; \real ) \to H_*(M\times M
; \real)$, and as $\Delta^*(a \otimes 1) = a$ we have that the
pullback $\Delta^*$ is surjective.  For $\tau \in \gr{S}_n$ the
map $f_\tau : (M^n)^\tau \to M^n$ is a composition of diagonal
maps, so we have that $f_{\tau *}$ is injective and that
$f_\tau^*$ is surjective. We can therefore apply proposition
\ref{proposition generators G} to the pair $(M^n, \gr{S}_n)$ to
get a set of generators. In what follows we will show that we can
reduce the set of generators by only considering the
transpositions.

\begin{lemma} \label{lemma generators symprod}
Let $\delta : M \to M^k$ be the diagonal inclusion and let
\begin{eqnarray*}
\sigma_{i}^k : M^{k-1} & \to & M^k\\
(x_1, \dots, x_{k-1}) & \mapsto & (x_1, \dots,x_{i-1}, x_i, x_i,
x_{i+1}, \dots x_{k-1})
\end{eqnarray*}
be the inclusion that repeats the $i$-th coordinate. Then in
cohomology
$$\delta_! 1 = \prod_{j=1}^{k-1} ({\sigma^k_{j}}_! 1).$$
\end{lemma}

\begin{proof}
We will proceed by induction on $k$. When $k=2$ the formula is
true because $\delta = \sigma_1^2$. Assume that we have shown the
formula for $k=n$ and let's try to show it for $k=n+1$. Consider
the following diagram of inclusions
$$\xymatrix{ M^n \ar[rr]^{\sigma_n^{n+1}} && M^{n+1} \\
M \ar[u]_\delta \ar[urr]_{\delta'} }$$ and using the properties of
the pushfoward we have,
\begin{eqnarray*}
\delta'_! 1 &=& {\sigma_n^{n+1}}_!(\delta_! 1) \\
& = & {\sigma_n^{n+1}}_!\left((\sigma_n^{n+1})^* ((\delta_! 1)
\otimes 1) \right)\\
& = & ({\sigma_n^{n+1}}_!1) ((\delta_!1) \otimes 1) \\
& = & ({\sigma_n^{n+1}}_!1) \left( \prod_{j=1}^n [({\sigma_j^n}_! 1)
\otimes 1] \right)\\
& = & ( {\sigma_n^{n+1}}_!1) \prod_{j=1}^{n-1} ({\sigma_j^{n+1}}_! 1)
\\
& = & \prod_{j=1}^{n} ({\sigma_j^{n+1}}_! 1).
\end{eqnarray*}
By the induction hypothesis the lemma follows.
\end{proof}

If we take the cycle $\alpha=(k, k-1, \dots, 2,1) \in \gr{S}_n$
and the transpositions $\tau_i = (i, i+1)$ we have that $f_{\alpha
!} 1 = \delta_! 1$ and $f_{\tau_i !} 1 = {\sigma_i^k}_!1$. By
lemma \ref{lemma generators symprod} we can conclude that for the
 cycle $\alpha$ of size $k$  which is the composition of $k-1$
transpositions $\tau_1 \dots \tau_{k-1}$ we have that
$$ f_{\alpha !} 1_\alpha = \prod_{i=1}^{k-1} (f_{\tau_i !} 1_{\tau_i}).$$

Therefore we can reduce the set of generators of the virtual
cohomology of the pair $(M^n, \gr{S}_n)$ by considering only the
transpositions.

\begin{proposition}
The ring $f(H^*_{virt}(M^n, \gr{S}_n ; \real))$ is generated by
the set
$$W= H^*(M^n ; \real)[1_{\gr{S}_n}] \cup \{(f_{\tau !}1_\tau) \tau \ | \ \tau
\in \gr{S}_n  \ \mbox{\rm{is a transposition}} \}$$ as a subring
of $H^*(M^n; \real)[\gr{S}_n]$.\\
If $H^*(M ; \integer)$ is torsion free, the same result holds but
with integer coefficients.
\end{proposition}

\begin{example}
Let's consider the pair $(M^n, \gr{S}_n)$ with $M=\complex P^m$.
We therefore have that $$H^*(M^n ; \integer) \cong \integer[x_1,
\dots , x_n]/\langle x_1^{m+1}, \dots , x_n^{m+1}\rangle$$ and we
only need to calculate the pushforward of the diagonal inclusion
$\Delta : M \to M\times M$.

By the Kunneth isomorphism we have that $$H_*( \complex P^m \times
\complex P^m ; \integer) \cong H_*(\complex P^m ; \integer)
\otimes H_*(\complex P^m ; \integer)$$ and $H^*(\complex P^m
\times \complex P^m; \integer) \cong \integer[x_1, x_2]/ \langle
x_1^{m+1}, x_2^{m+1} \rangle$. Let's take $H^*(\complex P^m ;
\integer) \cong \integer[y]/ \langle y^{m+1} \rangle$ and denote
by $[\complex P^i] \in H_{2i}(\complex P^m ;\integer)$ the
generator of the homology in degree $2i$ given by the inclusion
$\complex P^i \to \complex P^m$, $[z_0 : \dots : z_i] \mapsto [z_0
: \dots : z_i : 0 : \dots : 0]$.

By Poincar\'e duality (see \cite[page 213]{Hatcher}) we know that
the homology class $[\complex P^i]$ is dual to the cohomology
class $y^{m-i}$. Now,  $\Delta^* x_1 = \Delta^* x_2 = y$ implies
that $\Delta^* x_1^j x_2^{m-j} = y^m$, and as the classes
$\{x_1^jx_2^{m-j} | 0 \leq j \leq m \}$ generate the degree $2m$
cohomology of $\complex P^m \times \complex P^m$, we have by
Poincar\'e duality that $$ \Delta_* [\complex P^m] = \sum_{j=0}^m
[\complex P^j] \otimes [\complex P^{m-j}].$$Therefore we can
conclude that
 $\Delta_! 1 =\sum_{j=0}^m x_1^jx_2^{m-j}$.

 Then, for the transposition
$\tau=(k,l)$ we get that
$$f_{\tau !} 1 = \sum_{j=0}^{m} x_k^jx_l^{m-j}$$
and we can conclude that $H^*_{virt}((\complex P^m)^n, \gr{S}_n ,
\integer)$ is isomorphic to the subring of $\integer[x_1, \dots ,
x_n]/\langle x_1^{m+1}, \dots , x_n^{m+1}\rangle[\gr{S}_n]$
generated by
$$W=\{ (1 \ 1_{\gr{S}_n}) \} \cup \{(x_i \ 1_{\gr{S}_n}) \ | \ 1 \leq i \leq n \} \cup
\left\{\left(\sum_{j=0}^m x_k^j x_l^{m-j}\  (k,l)\right)  \ | \ 1
\leq k < l \leq n \right\}$$
\end{example}

\begin{example}
Let's pay particular attention to the case on which $n=2$ and $M$
is a connected, differentiable, compact and closed manifold of
dimension $d$. Denote by $\Omega \in H^d(M ; \real)$ the generator
of the top cohomology, then we have that
$$(\Delta_!1)(\Delta_!1) = \chi(M) \Omega \otimes \Omega$$
where $\chi(M)$ is the Euler number of $M$ (see \cite[Pro.
11.24]{BottTu}). Moreover, by the properties of the pushforward in
cohomology we have that for any $\alpha, \beta \in H^*(M ; \real)$
$$(\Delta_!1) (\alpha \otimes \beta) = \Delta_!(\Delta^*(\alpha
\otimes \beta)) = \Delta_!(\alpha\beta),$$ hence, if
$\alpha_1\beta_1 = \alpha_2\beta_2$ we have that
$$(\Delta_!1)(\alpha_1 \otimes \beta_1) = (\Delta_!1)(\alpha_2
\otimes \beta_2)$$
and if $deg(\alpha) + deg(\beta) > d$
$$(\Delta_!1) (\alpha \otimes \beta) =0.$$

 If we consider the ring $$\ H^*(M ; \real)^{\otimes 2}
 [u]$$ where $u$ represents the element $(\Delta_!1)$ then
we can se that $u^3=0$ and $u^2 - \chi(M) \Omega \otimes \Omega$. The annihilator ideal of
$u$ is generated by the elements $\alpha \otimes \beta $ where $deg(\alpha) + deg(\beta) > d$, and
by the elements $(\alpha_1 \otimes \beta_1 - \alpha_2 \otimes \beta_2)$ where $\alpha_1\beta_1 =
\alpha_2\beta_2$.

In the case that $M=\complex P^m$ we can see that
$H_{virt}^*((\complex P^m)^2, \gr{S}_2 ; \integer)$ is isomorphic
to
$$\integer[x,y,u]/\langle x^{m+1}, y^{m+1}, u^2 - (m+1)x^my^m,
u(x-y) \rangle$$ where $u^3=0$ because $u^3 = (m+1)x^my^mu = (m+1)
x^{m-1} y^{m+1} u =0$.

In the case that $m=1$ we have that $H_{virt}^*((\complex P^1)^2,
\gr{S}_2 ; \integer)$ is
$$\integer[x,y,u]/\langle x^{2}, y^{2},  u^2 - 2xy,
u(x-y) \rangle. $$ The $\integer/2$ invariants are generated as an
$\real$-module by $x+y, u, xy$ and $ux$. Therefore, if we take $w
= x +y$, then $w^2 = 2xy= u^2$, $2ux=uw$ and $w^3=0$. So  we can
conclude  that
$$H^*_{virt}([(\complex P^1)^2/ \gr{S}_2]; \real) \cong \real[w,u]/\langle w^{3}, u^3, u^2 - w^2
 \rangle .$$

\end{example}




\bibliographystyle{alpha}
\bibliography{VirtualCohomology}
\end{document}